\documentclass[a4paper,10pt]{article}

\usepackage{graphicx}      

\usepackage{graphicx} 

\usepackage{savesym}
\savesymbol{AND}
\usepackage{algorithm, algorithmic, setspace}
\usepackage{graphicx} 
\usepackage{amsmath} 
\usepackage{amssymb}  
\usepackage{amsfonts}
\usepackage{color}
\usepackage[english]{babel}

\newenvironment{proof}{\paragraph{Proof}}{\hfill$\square$}

\newcommand{\argmin}[1]{\underset{#1}{\mathrm{argmin\,}}}
\newcommand{\xmin}{x^{\star}}
\newcommand{\Xmin}{X^{\star}}
\newcommand{\Zmin}{Z^{\star}}
\newcommand{\xtrue}{\widetilde{x}}
\newcommand{\ztrue}{\widetilde{z}}

\newcommand{\R}{\mathbb{R}}
\newcommand{\Syp}{\mathbb{S}_+}
\newcommand{\Sy}{\mathbb{S}}

\newcommand{\one}{\mathbf{1}}

\newcommand{\tr}{\mathrm{tr}}

\newtheorem{proposition}{Proposition}

\newtheorem{remark}{\noindent \textbf{Remark}}{\normalfont}{\normalfont}
{\normalfont}{\normalfont}
 {\normalfont }{\normalfont}
{\normalfont}{\normalfont}

\begin{document}

\title{Enhancing low-rank solutions in semidefinite relaxations of Boolean quadratic problems} 

\author{V. Cerone, S. M. Fosson$^{*}$, D. Regruto}
\maketitle
\begin{abstract}  
$~$ Boolean quadratic optimization problems occur in a number of applications. Their mixed integer-continuous nature is challenging, since it is inherently NP-hard. For this motivation,  semidefinite programming relaxations (SDR's) are proposed in the literature to approximate the solution, which recasts the problem into convex optimization. Nevertheless, SDR's do not guarantee the extraction of the correct binary minimizer.  In this paper, we present a novel approach to enhance the binary solution recovery. The key of the proposed method is the exploitation of known information on the eigenvalues of the desired  solution. As the proposed approach yields a non-convex program, we develop and analyze an iterative descent strategy, whose practical effectiveness is shown via numerical results.
\end{abstract}



\section{Introduction}
Boolean quadratic optimization refers to those minimization problems with quadratic cost functional and binary variables. It is a long-time problem with a number of applications in different scientific areas. To mention some examples, it is encountered in maximum cut (MAX-CUT) problems in graphs, see, \cite{pal12}, knapsack problems, see \cite{kel04}, hybrid model predictive control, see \cite{nai17,bem18}, sensor selection, see \cite{she14,car13s}, medical imaging, see \cite{sch05, web05},  and binary compressed sensing, see \cite{fox19spl}.

Boolean quadratic optimization is challenging as it is NP-hard, even when the cost functional is convex, due to the integer nature of the variables. To overcome this drawback, different relaxations are proposed in the literature that approximate the correct solution. The most popular approach consists in the semidefinite programming  relaxation (SDR), known also as Shor relaxation, introduced by \cite{sho87} and \cite{lov90}. In a nutshell, the key idea at the basis of SDR is to embed the variable matrix $xx^T$, $x\in\R^n$, into the space of $n\times n$ symmetric positive semidefinite  matrices. The rationale is that a complete description of the cone of symmetric positive semidefinite matrices is available, while this does not hold for the convex polytope of matrices $xx^T$, as illustrated, e.g., by \cite{nes98}. The tightness of SDR is analyzed, e.g., in \cite{goe95,nes98}.

More recently, a substantial evolution on SDR has been provided by the Lasserre's theory, which states that the global optimum of a polynomial optimization problem can be achieved by solving a hierarchy of SDR's; we refer the reader to  \cite{las01,lasbook} for details, and to \cite{las00} for the specific case of quadratic functionals.
Theorem 4.2 in \cite{las01} and further elaborations in Chapter 6 of \cite{lasbook} provide sufficient conditions to get the global minimum and a global minimizer of the original problem from an  SDR of a certain relaxation order. As to the Shor relaxation (or equivalently, SDR of order 1) of a Boolean quadratic problem, if the SDR solution is 1-rank, then the global minimum is achieved and the desired binary minimizer can be extracted, as discussed in \cite{fox19spl}. For this motivation, methods to minimize or reduce the rank in semidefinite programming are relevant for SDR; we refer to \cite{lem16} for a complete overview on the topic. Since the rank minimization is an NP-hard problem, the minimization of the nuclear norm is often exploited, which provides the tightest convex relaxation of the problem, as proven by \cite{rec10}. As to symmetric positive semidefinite matrices, the nuclear norm is equal to the trace, and in turn to the sum of the eigenvalues: then, the minimization of the nuclear norm has a sparsifying effect on the eigenvalues vector, as the $\ell_1$-norm is  a suitable convex relaxation of the sparsity a vector, see, e.g., \cite{tib96}. Then, a nuclear norm penalty encourages a low-rank solution. Nevertheless, this method does not apply to some Boolean quadratic problems, e.g., MAX-CUT, see \cite{las16}, where the trace is constant. For these problems, an effective method is the  $\log$-$\det$ heuristic, proposed by \cite{faz03}.

In this paper, we propose a novel method to enhance SDR for Boolean quadratic optimization,  by supporting 1-rank solutions. The key idea is to exploit the knowledge of the eigenvalues of the desired solution are known (in particular, only one eigenvalue is non-null). Then, we promote  1-rank solutions by maximizing the energy of  the eigenvalues vector, which we prove to have a sparsifying effect on the eigenvalues, when the trace is constant. We illustrate the proposed approach for two SDR's, with slightly different spectral properties. Furthermore, as the method yields a non-linear and non-convex program, a suitable iterative descent algorithm is developed to search the global minimum. The overall strategy is analyzed, and its effectiveness is illustrated through numerical simulations.

The paper is organized as follows. In Sec. \ref{sec:PS}, we formally state the problem. In Sec. \ref{sec:PA}, we present and analyze the proposed approach. In Sec. \ref{sec:NR}, we show some numerical results; finally, we draw some conclusions in Sec. \ref{sec:C}.
\section{Problem Statement}\label{sec:PS}
Let $\Sy^n$ be the space of $n\times n$ symmetric matrices, and $\Syp^n\subset \Sy_n$ the subspace of positive semidefinite matrices. 
We consider Boolean quadratic problems of the kind 
\begin{equation}\label{prob:1}
\min_{x\in\{0,1\}^n} x^T C x + 2d^Tx {\text{ s. t. }} Ax=b
\end{equation}
where $C\in\Sy^n$, $d\in\R^n$, $A\in\R^{m,n}$, and $b\in\R^m$.
%
As in \cite{las16}, we exploit the Lagrangian formulation
\begin{equation}\label{prob:2}
\min_{x\in\{0,1\}^n} x^T C x + 2d^Tx+\mu\|Ax-b\|_2^2
\end{equation}
where $\mu>0$. If $\mu$ is sufficiently large, problems \eqref{prob:1} and \eqref{prop:2} are equivalent, see \cite{las16}. 
By adding a slack variable $x_0=1$  and by redefining $x=(1,x_1,\dots, x_n)^T$, we rewrite \eqref{prob:2} as the augmented problem
\begin{equation}\label{prob:3}
\min_{x\in\{0,1\}^{n+1}} x^T Q x
\end{equation}
where $Q = \left(\begin{array}{c|c}
           0&d^T\\
           \hline
           d&C\\
          \end{array}\right)+\mu(-b| A)^T(-b | A).$
The solution to this NP-hard problem can be searched by solving the associated SDR, see \cite{nes98} for details:
\begin{equation}\label{sdr:1}
\begin{split}
\min_{X\in\Syp^{n+1}} &\langle Q,X \rangle\\
\text{s. t. }&X_{0,0}=1; ~~~~X_{i,i}=X_{0,i},~~i=1,\dots,n.
\end{split}
\end{equation}
The constraint $X_{i,i}=X_{0,i}$ represents the constraint $x_i=x_i^2$, which holds for any $x_i\in\{0,1\}$. 

An alternative, MAX-CUT approach is studied in \cite{las16}: given $\one=(1,\dots,1)^T\in\R^{n+1}$, by substituting  $z=2x-\one\in\{-1,1\}^{n+1}$  in \eqref{prob:3}, we obtain $x^TQx=\frac{1}{4} (z+\one)^TQ(z+\one) $. Then, problem \eqref{prob:3} is equivalent to 
\begin{equation}\label{prob:4}
\min_{z\in\{-1,1\}^{n+1}} z^T R z
\end{equation}
where $R = Q+ \left(\begin{array}{c|c}
           0&\one^T Q\\
           \hline
           Q\one&\mathbf{0}\\
          \end{array}\right)$.
Problem \eqref{prob:4} is a MAX-CUT problem, and the associated SDR is
\begin{equation}\label{sdr:2}
\begin{split}
\min_{Z\in\Syp^{n+1}} &\langle R ,Z \rangle\\
&\text{s. t. } Z_{i,i}=1,~~~i=0,\dots, n.
\end{split}
\end{equation}
 The constraint $Z_{i,i}=1$ represents the fact that $z_i^2=1$ whenever $z_i\in\{-1,1\}$. 
 
The aim of this paper is the recovery of the correct binary minimizers  of problems \eqref{prob:3} and \eqref{prob:4} by starting from their SDR's \eqref{sdr:1} and \eqref{sdr:2}. As described, e.g, in \cite{las16} and \cite{lem16}, the correct binary minimizers can be extracted if and only if the SDR solutions are 1-rank. Therefore, we develop a strategy to reduce, possibly minimize, the rank, by exploiting specific features of the Boolean setting. 
\section{Concave penalization exploiting information on the eigenvalues}\label{sec:PA}
In this section, we develop the proposed strategy to promote 1-rank solutions to SDR's \eqref{sdr:1} and \eqref{sdr:2}. Specifically, we propose a suitable cost functional, and we illustrate the algorithms used to minimize it.
\subsection{Proposed cost functional}
Let us consider SDR \eqref{sdr:2}, and let us call $\Zmin\in\Syp^{n+1}$ the desired 1-rank solution. The rank of $\Zmin$ corresponds to the number of non-null eigenvalues; thus, the  eigenvalues vector of $\Zmin$, denoted as $v$, is 1-sparse, i.e., $v=(v_0,0,\dots,0)$, $v_0> 0$. We also remark that  $\Zmin$ necessarily has all the components in $\{-1,1\}$, and $\frac{1}{2}((\Zmin_{0,1},\dots, \Zmin_{0,n})^T+\one)$ is the exact minimizer of \eqref{prob:2}.  \cite{rec10,lem16} show that a sparsifying effect on $v$ can be obtained by penalizing the trace of the variable $Z$, hereafter denoted as $\tr(Z)$, which corresponds to the $\ell_1$-norm of the eigenvalues, and, in turn, the $\ell_1$-norm is a suitable convex relaxation of the sparsity of a vector. Nevertheless, this approach is not beneficial for SDR \eqref{sdr:2}, as $\tr(Z)=n+1$ by construction; therefore, it makes no sense to penalize a constant quantity. 

However, we observe that the information $\tr(\Zmin)=\sum_{i=0}^{n} v_i=n+1$ can be exploited to state that $v_0=n+1$. In other terms, not only we know that $v$ is 1-sparse, but also that its components belong to the binary alphabet $\{0,n+1\}$. Then, we wonder how to force the solution to have $v=\{n+1,0,\dots,0\}$.

The key idea is as follows. Given a vector $v\in [0,n+1]^{n+1}$ with $\sum_{i=0}^n v_i=n+1$, we can force 1-sparsity by maximizing its energy $\|v\|_2^2$. This is straightforward to check in the two-dimensional case: let us consider $(v_0,v_1)\in[0,2]^2$ with $v_0+v_1=2$;  the maximum of $v_0^2+v_1^2=v_0^2+(2-v_0)^2$ is achieved at the boundaries, that is, at $(0,2)$ or $(2,0)$.  This reasoning can be extended to any dimension.

Given this principle, we search a method to maximize the energy of $v$ within SDR \eqref{sdr:2}. We notice that 
$\|v\|_2^2=\tr(ZZ)=\langle Z, Z \rangle.$
Then, we propose to add  a term that penalizes $-\langle Z, Z \rangle$ in the cost functional of \eqref{sdr:2}, i.e.,
\begin{equation}\label{sdr:2c}
\begin{split}
\min_{Z\in\Syp^{n+1}} &\langle R ,Z \rangle - \lambda\langle Z,Z\rangle\\
&\text{s. t. } Z_{i,i}=1,~~~i=0,\dots, n\\
\end{split}
\end{equation}
where $\lambda>0$ is a design parameter that can be assessed by cross-validation. Interestingly, if a global minimizer of \eqref{sdr:2} is binary, then, for any $\lambda>0$, the global minimizer of \eqref{sdr:2c} is exact, as illustrated in the following proposition.
%

\begin{proposition}\label{prop:1}
Let $\ztrue\in\{-1,1\}^{n+1}$ be the correct solution to problem \eqref{prob:4}. Let us assume that SDR \eqref{sdr:2} has the desired binary, 1-rank solution $\ztrue\ztrue^T$ among its global minimizers. Then, the minimizer of \eqref{sdr:2c} is $\Zmin=\ztrue\ztrue^T$, for any $\lambda>0$. Moreover, $\ztrue=(1, \Zmin_{1,0},\dots, \Zmin_{n,0})^T$. 
\end{proposition}
\begin{proof}
Since $\tr(Z)=n+1$, then $\langle Z,Z\rangle\leq (n+1)^2$. In particular, the maximum $(n+1)^2$ is achieved if and only if the eigenvalues of $Z$ are $v=(n+1,0,\dots,0)$. In this case, $Z$ is 1-rank, and necessarily $Z=\ztrue\ztrue^T$.
Moreover, since by assumption $\ztrue\ztrue^T$ minimizes $\langle R ,Z \rangle$, then it is the unique global minimizer of \eqref{sdr:2c}.
\end{proof}

A similar approach can be applied to SDR \eqref{sdr:1}, with a considerable difference: $\tr(X)$ is not priorly set by construction. Then, we proceed by considering two possible settings. 

In the first setting, we assume to know the sparsity level $k$ of the true solution of problem \eqref{prob:3}. This implies that the desired solution to SDR \eqref{sdr:1} has trace equal to $(k+1)^2$ and eigenvalues $v=(k+1,0,\dots,0)$. The sparsity level is known in many applications, e.g., in compressed sensing, see  \cite[Sec. III]{fou13}, or in sensor selection, where $k$ is the number of used sensors, see \cite{she14}. In other cases, an unknown $k$ can be estimated through ad hoc techniques, see \cite{rav18}. Given $k$, we propose to modify SDR \eqref{sdr:1} as follows to leverage the binary nature of the eigenvalues:
\begin{equation}\label{sdr:1c}
\begin{split}
\min_{X\in\Syp^{n}} &\langle Q,X \rangle+\lambda\left[ (k+1)\tr(X)-\langle X,X\rangle\right]\\
\text{s. t. }&X_{0,0}=1;~~~~X_{i,i}=X_{0,i},~~i=1,\dots,n.
\end{split}
\end{equation}
%

\begin{proposition}\label{prop:2}
Let $\xtrue\in\{0,1\}^{n+1}$ be the correct solution to problem \eqref{prob:3}. Let us assume that SDR \eqref{sdr:1} has the desired binary, 1-rank solution $\xtrue\xtrue^T$ among its global minimizers. Then, the minimizer of problem \eqref{sdr:1c} is $\Xmin=\xtrue\xtrue^T$, for any $\lambda>0$. Moreover, $\xtrue=(1,    \Xmin_{1,0},\dots, \Xmin_{n,0})^T$. 
\end{proposition}
\begin{proof}
If $v$ is the eigenvalues vector of $X$, we have
\begin{equation}
 \tr(X)-\langle X,X\rangle = \sum_{i=1}^n [(k+1) v_i-v_i^2]\geq 0
\end{equation}
since $v_i\in [0,k+1]$ for each $i=1,\dots,n$.
Moreover, $\tr(X)-\langle X,X\rangle = 0 ~~\Longleftrightarrow~~ v_i\in \{0,k+1\}$ 
for each  $i=1,\dots,n$. Then, we can conclude that the global minimum is achieved for $v\in\{0,k+1\}^{n+1}$. Eventually, since $\sum_{i=1}^n v_i =k+1$, we have $v=(k+1,0,\dots,0)$, which implies the 1-rank solution, and $\Xmin=\xtrue\xtrue^T$.
\end{proof}

In the second setting, we assume that $k$ is unknown. In this case, we propose to replace $k$ by $n\geq k$. This approximate procedure takes a larger weight on the term $\tr(X)$, which might be advantageous when $k\ll n$: in fact, in this setting, not only $\tr(\Xmin)=\sum_{i=1}^{n+1} v_i$, but also $\tr(\Xmin)=1+\sum_{i=1}^n \xmin_i$, where $\xmin\in [0,1]^n$ is the final estimation of $\xtrue$, extracted from the diagonal of $\Xmin$: $\xmin=(\Xmin_{1,1},\dots, \Xmin_{n,n})^T$. In conclusion, by penalizing $\tr(X)$ we obtain a sparsifying effect both on the eigenvalues and on the solution; for this motivation, we expect better performance when $k\ll n$. 

\begin{remark}\label{rem:1}
The considered problem bears some similarities with phase retrieval from  Fourier transform magnitude, see, e.g., \cite{jag13,jag17}. An effective approach to phase retrieval is the embedding of the unknown vector $x$ into a higher dimensional space by the transformation $X=xx^T$. This approach, called lifting, shares the same principle of the SDR's illustrated in this paper. Similarly to the MAX-CUT problem, the penalization of the trace is not effective for phase retrieval, as the energy of $x$ is fixed, see \cite[Sec. III]{jag13}. In \cite{jag13}, this issue is overcome by penalizing the term $\log\det(X+\epsilon I)$, where $\epsilon>0$ is a small design parameter necessary for boundedness. This term is a concave surrogate of the rank, whose effectiveness is discussed in \cite{faz03}. Our approach  consists in a  concave penalization as well. However, differently from the $\log$-$\det$ heuristic, it  exploits the known binary eigenvalues; moreover, its practical implementation is computationally less complex, as illustrated in Remark \ref{rem:2}.
\end{remark}
\subsection{Descent algorithms}
Problems \eqref{sdr:1c} and \eqref{sdr:2c} are well posed, that is, their unique global minima are the correct binary solutions. However, they introduce the concave term $-\langle X, X \rangle$, which makes the problem non-linear and non-convex. For this motivation, we propose an iterative descent algorithm to search the minimum, which, although sub-optimal, is effective in practice. The idea is to replace the concave term $\langle X,X\rangle$ with $\langle G,X\rangle$, where $G$ is an available estimate of $X$. By assuming $G$ fixed, the penalty is linear in $X$, then the whole cost functional is linear. Then, we propose an iterative procedure: we start from an initial estimate of the solution, denoted by $X_0$, which can be assessed, for example, by solving the non-penalized problems \eqref{sdr:1} and \eqref{sdr:2}; in turn, we solve the penalized problems until convergence is reached. The overall procedure is summarized in algorithms \ref{alg:1} and \ref{alg:2} for problems \eqref{sdr:1c} and \eqref{sdr:2c}, respectively . In Algorithm \ref{alg:1}, $h$ is  equal to $k+1$ when $k$ is known, and $n+1$ otherwise.

\begin{algorithm}
\setstretch{1.2}
     \renewcommand{\algorithmicrequire}{\textbf{Input:}}
    \renewcommand{\algorithmicensure}{\textbf{Output:}}
	\caption{Descent algorithm for Problem \eqref{sdr:1c} }\label{alg:1}
	\begin{algorithmic}[1]
		\REQUIRE $Q, \lambda>0$;
		\STATE $X_0=$ solution of \eqref{sdr:1}
		\FORALL{$t=1,\dots,T$}
		\STATE $X_{t}=\argmin{X\in\Syp^{n+1}} \langle Q,X \rangle  +\lambda\left[h\tr(X)- \langle X_{t-1}, X\rangle\right]$,$~~$ s. t. $X_{0,0}=1$, $~X_{i,i}=X_{0,i},~i=1,\dots,n$
		\ENDFOR
	\end{algorithmic}
\end{algorithm}

\begin{algorithm}
\setstretch{1.2}
     \renewcommand{\algorithmicrequire}{\textbf{Input:}}
    \renewcommand{\algorithmicensure}{\textbf{Output:}}
	\caption{Descent algorithm for Problem \eqref{sdr:2c}  }\label{alg:2}
	\begin{algorithmic}[1] 
		\REQUIRE Input: $R, \lambda>0$;
		\STATE $X_0=$ solution of \eqref{sdr:2}
		\FORALL{$t=1,\dots,T$}
		\STATE $X_{t}=\argmin{X\in\Syp^{n+1}} \langle R,X \rangle-\lambda \langle X_{t-1}, X\rangle$, s. t.  $X_{i,i}=1$, $i=0,\dots,n$
		\ENDFOR
	\end{algorithmic}
\end{algorithm}
In this way, we solve a sequence of semidefinite programming problems. This is not guaranteed to get the global minimum, while it is guaranteed to provide a non-increasing cost functional sequence.
\vskip0.3cm

\begin{proposition}
Let us define $f(X):=\langle Q,X\rangle +\lambda h\tr(X)$ from \eqref{sdr:1c} (respectively, $f(X):= \langle R,X\rangle$ in \eqref{sdr:2c}), and  $F(X_t):=f(X_t)-\frac{\lambda}{2}\langle X_t,X_t \rangle$. 

By applying Algorithm \ref{alg:1} (respectively, Algorithm \ref{alg:2}), $F(X_t)$ is a non-increasing function.
\end{proposition}
\begin{proof}
 Since $X_t$ is the minimizer, 
 \begin{equation}
\begin{split}
f(X_t)-\lambda\langle X_{t-1}, X_{t}\rangle &\leq f(X_{t-1})-\lambda\langle X_{t-1}, X_{t-1}\rangle.
\end{split}
\end{equation}
On the other hand, if $A,B\in\Syp^n$, then  $2\langle A,B \rangle \leq \langle A,A \rangle + \langle B,B \rangle$, see, e.g., \cite{zho14}. Therefore, for any $t$,
$2\langle X_{t-1},X_t \rangle - \langle X_{t-1},X_{t-1} \rangle \leq   \langle X_t,X_t \rangle$, which implies 
\begin{equation*}
\begin{split}
F(X_t)& \leq f(X_t)-\lambda\langle X_{t-1},X_t \rangle + \frac{\lambda}{2}\langle X_{t-1},X_{t-1}\rangle\\
& \leq f(X_{t-1})-\lambda \langle X_{t-1},X_{t-1}\rangle  +\frac{\lambda}{2}\langle X_{t-1},X_{t-1}\rangle\\
& \leq F(X_{t-1})-\frac{\lambda}{2} \langle X_{t-1},X_{t-1}\rangle  +\frac{\lambda}{2}\langle X_{t-1},X_{t-1}\rangle\\
&=F(X_{t-1}).
\end{split}
\end{equation*}
\end{proof}

Furthermore, if the proposed approach achieves a binary solution, by uniqueness it is guaranteed that this is exactly the desired global minimum. Conversely, when the obtained solution is not 1-rank, the final solution is not binary and it is guaranteed that it is not the correct solution. This awareness about achieving the correct solution provides the possibility of running again the algorithm by suitably changing the initialization, which may yield a better solution. 

\begin{remark}\label{rem:2}
The proposed descent strategy can be interpreted as a reweighting algorithm. In \cite{jag13,faz03}, the reweighting algorithm is derived as iterative minimization of the local linearization of the concave $\log$-$\det$ term. The convergence of the  $\log$-$\det$ method to a local minimum is discussed in \cite{faz03}. We mention that the proposed descent algorithm can be described under the reweighting viewpoint as well, and results from \cite{faz03,fox18spl} might be leverage to rigorously prove its convergence. The convergence analysis and a complete comparison to  $\log$-$\det$ heuristic are beyond the scope of this paper; however, some numerical comparisons are proposed in Sec. \ref{sec:NR}. 
\end{remark}
\section{Numerical results}\label{sec:NR}
In this section, we illustrate some numerical results, that support the effectiveness of the  proposed approach. We consider the following problem, as presented in \cite{fox19spl}: we aim to solve the underdetermined system 
\begin{equation}\label{m1}
b=A x,~~~~ x\in\{0,1\}^n,~~b\in\R^m,~A\in\R^{m,n},~~m< n
\end{equation}
under the assumption that the solution is unique in $\{0,1\}^n$. This linear problem is encountered in a number of applications, ranging from compressed sensing, see \cite{fli18,fox18asi,fox19} to tomography, see \cite{sch05,web05}.

In \cite{fox19spl}, problem \eqref{m1} is tackled in case of sparse $x$, which recasts into binary compressed sensing, and a Shor SDR is proposed to solve it. In previous works, other methods for binary compressed sensing are proposed, namely, relaxation over the convex hull in \cite{fli18}, $\ell_1$-reweighting algorithms in \cite{fox18asi}, alternating direction method of multipliers in \cite{fox19}, and difference of convex functions in \cite{sch05,web05}. SDR is shown to achieve better accuracy in \cite{fox19spl}.

Here, we recast problem \eqref{m1}  into  Boolean quadratic optimization by replacing the constraint $Ax=b$ with the the cost functional $\|Ax-b\|_2^2$. As we know that $\xtrue\in\{0,1\}^n$ is solution of $Ax=b$, the global minimum is null. 
Given, the cost functional $\|A x -b\|_2^2$, we homogenize it by adding the slack variable $x_0=1$, so that $\|A x -x_0b\|_2^2=x^TQx$ with $Q=(-b| A)^T(-b|A)$, and $x=(x_0,x_1,\dots, x_n)^T$. Then, we can apply the approach developed in Sec. \ref{sec:PA}.

Specifically, we compare the recovery accuracy of the proposed ``known binary eigenvalues'' (KBE) approach to the following known methods: SDR \eqref{sdr:1} and SDR \eqref{sdr:2}, nuclear norm heuristic, and $\log$-$\det$ heuristic, see \cite{faz03,lem16}.
The implemented nuclear norm and $\log$-$\det$ algorithms read as follows, respectively:
\begin{equation}\label{sdr:1b}
\begin{split}
\min_{X\in\Syp^{n+1}} &\langle Q,X \rangle+\lambda \tr(X)\\
&\text{s. t. }X_{0,0}=1,~~X_{i,i}=X_{0,i},~~i=1,\dots,n
\end{split}
\end{equation}
and, for $t=1,\dots,T$,
\begin{equation}\label{sdr:ld}
\begin{split}
&~~~~~X_{t}=\argmin{X\in\Syp^{n+1}} \langle Q,X \rangle  +\lambda\langle (X_{t-1}+\epsilon I)^{-1},X \rangle\\
&~~~~~\text{s. t. }X_{0,0}=1,~~X_{i,i}=X_{0,i},~~i=1,\dots,n.
\end{split}
\end{equation}
\begin{figure}[ht]
\begin{center}
\includegraphics[width=0.6\columnwidth]{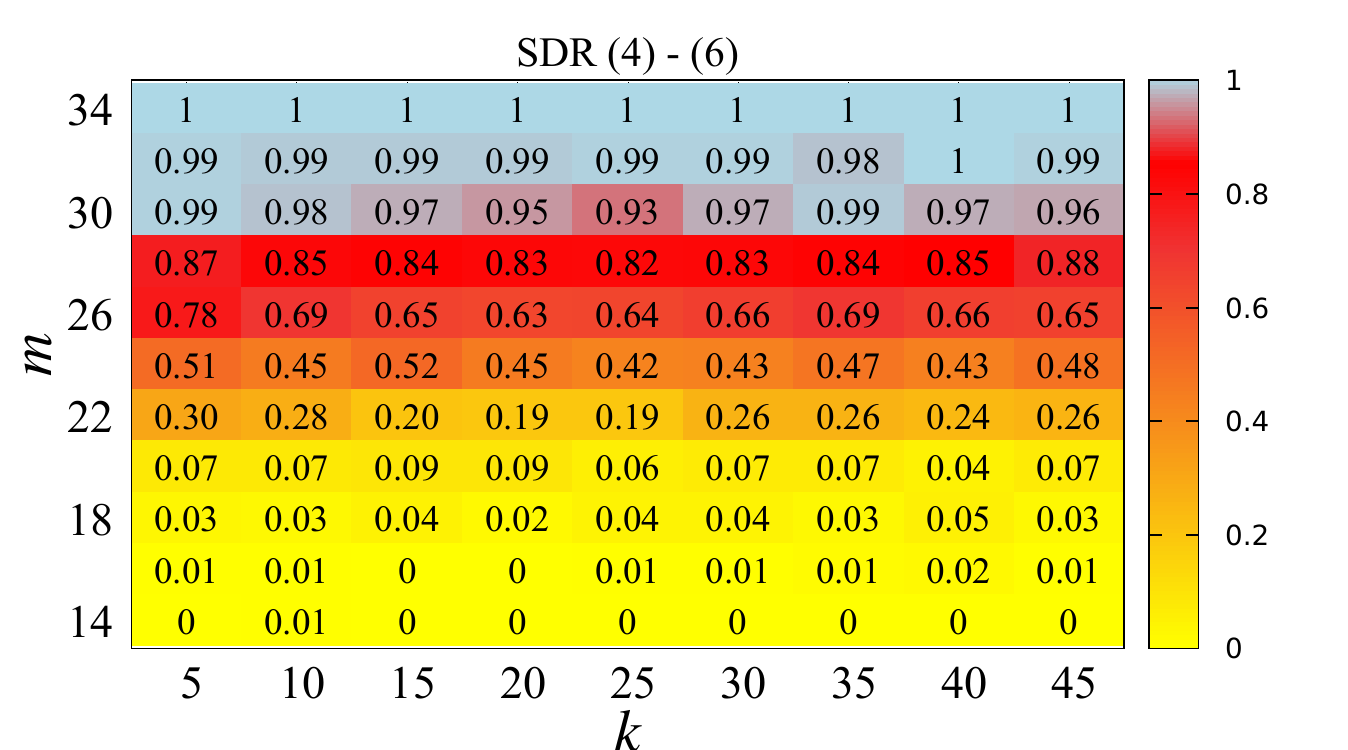}\\ 
\includegraphics[width=0.6\columnwidth]{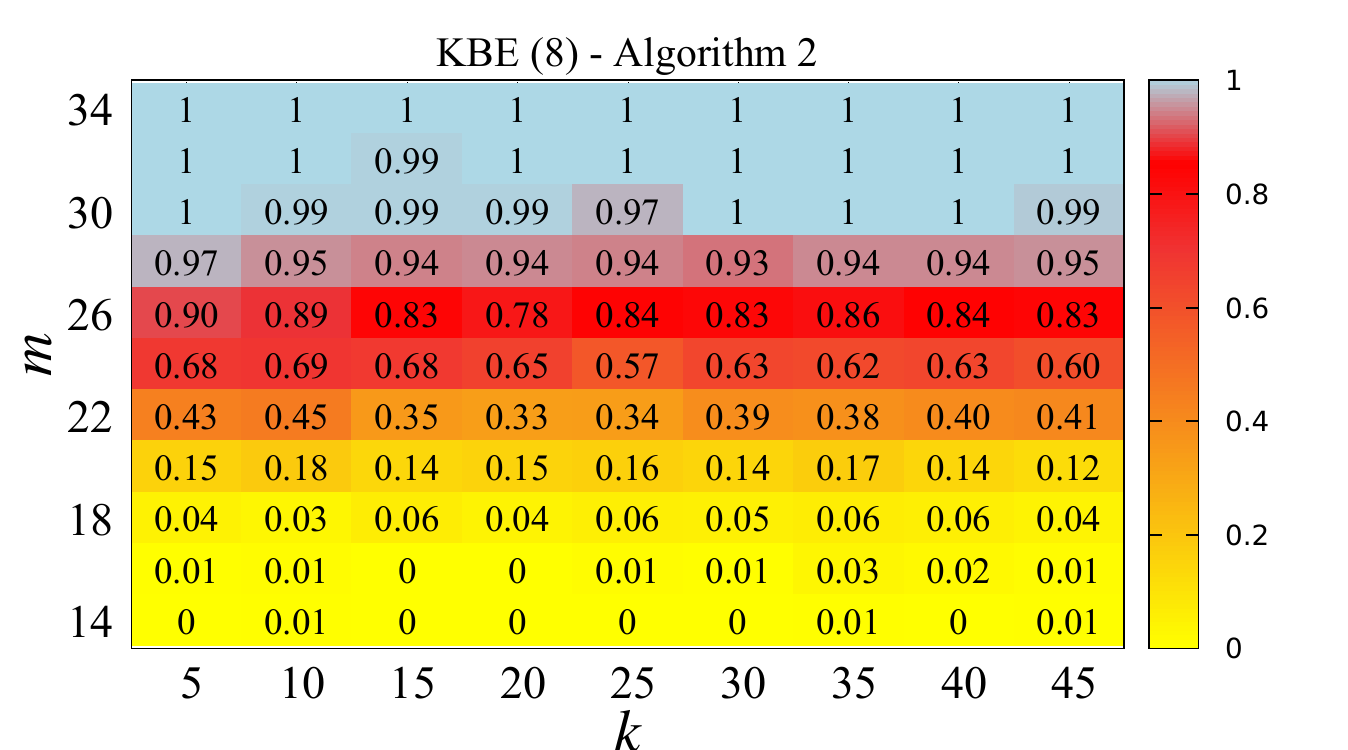}\\
\includegraphics[width=0.6\columnwidth]{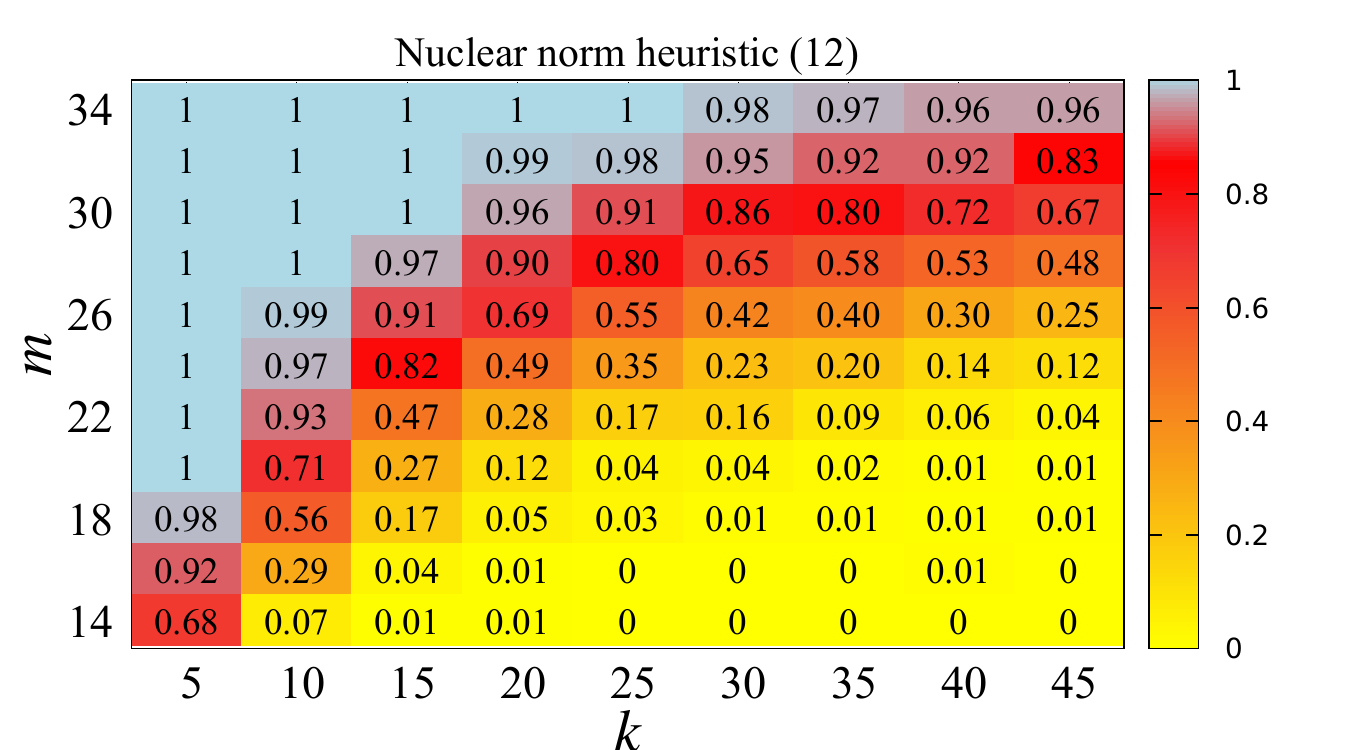}\\
\includegraphics[width=0.6\columnwidth]{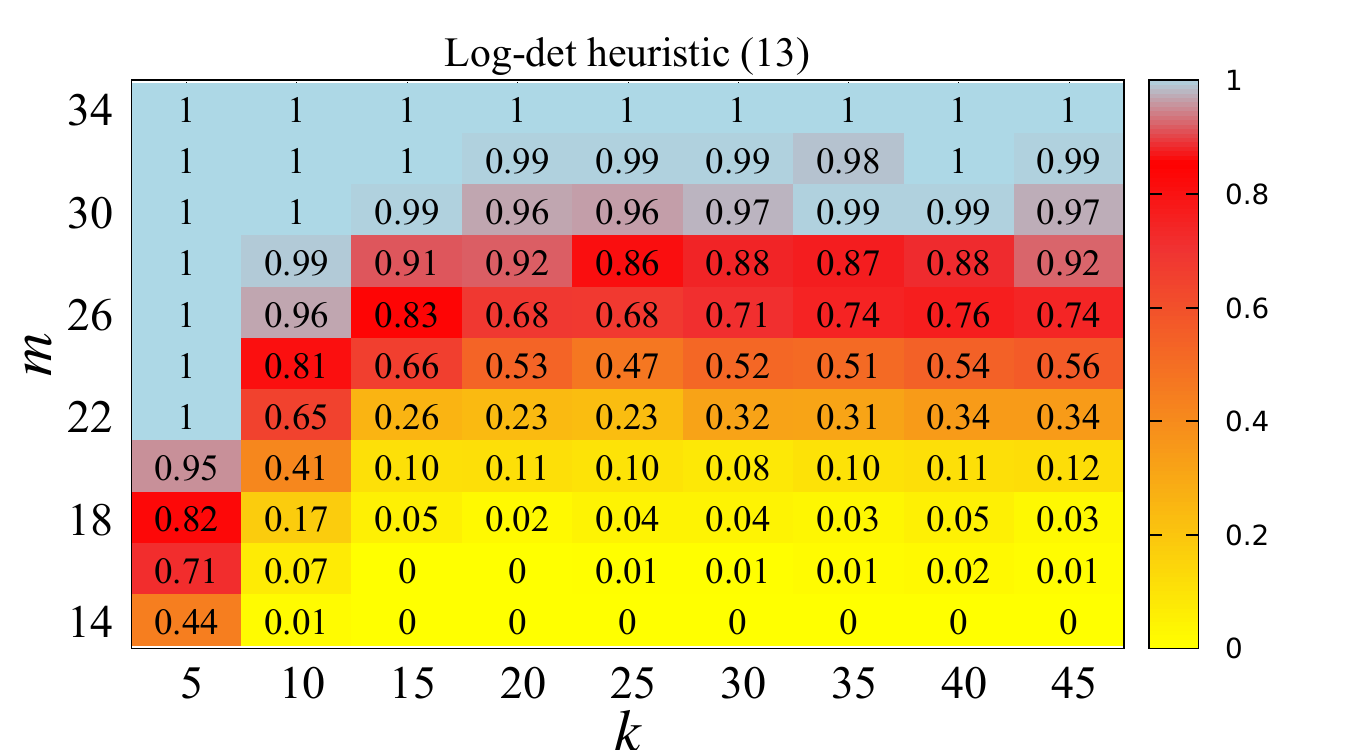}\\
\includegraphics[width=0.6\columnwidth]{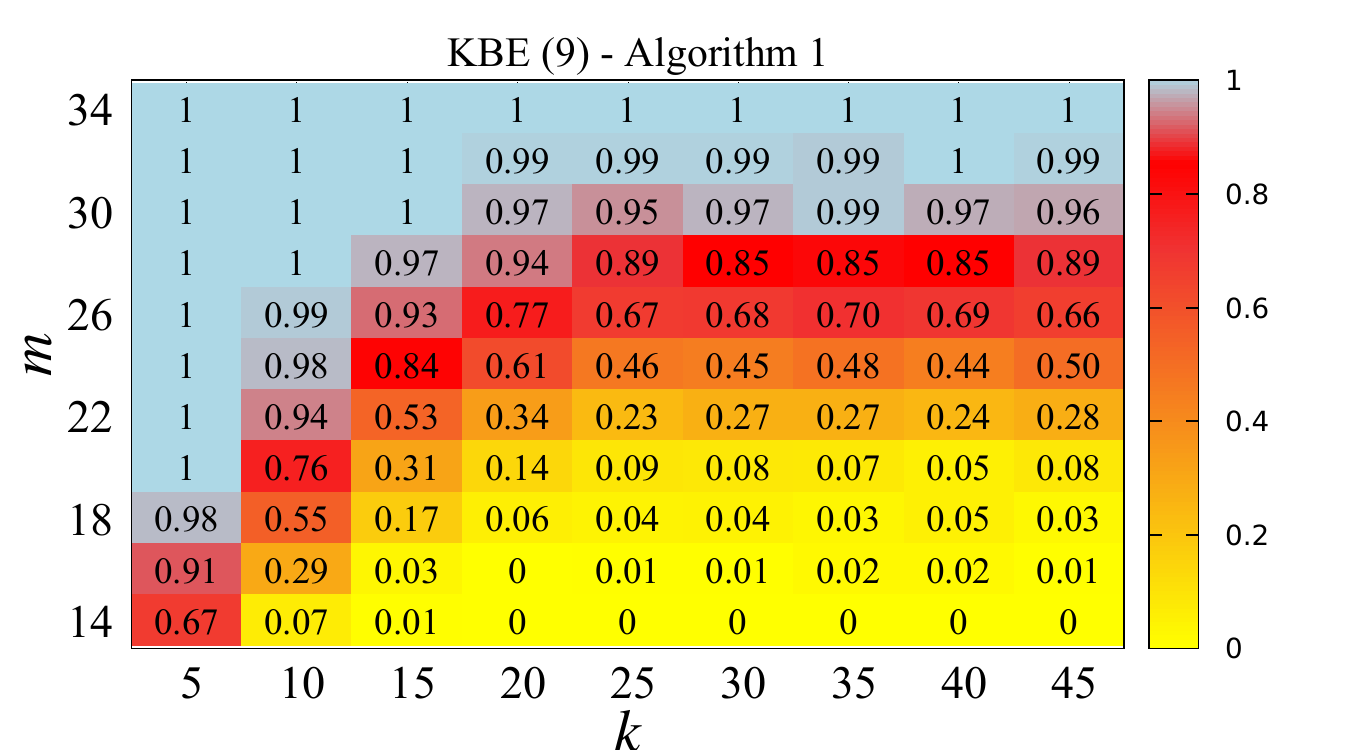}
\caption{Experiment 1 (unknown $k$): exact recovery rate.}
\label{fig:1}
\end{center}
\end{figure}
\begin{figure}[ht]
\begin{center}
\includegraphics[width=0.6\columnwidth]{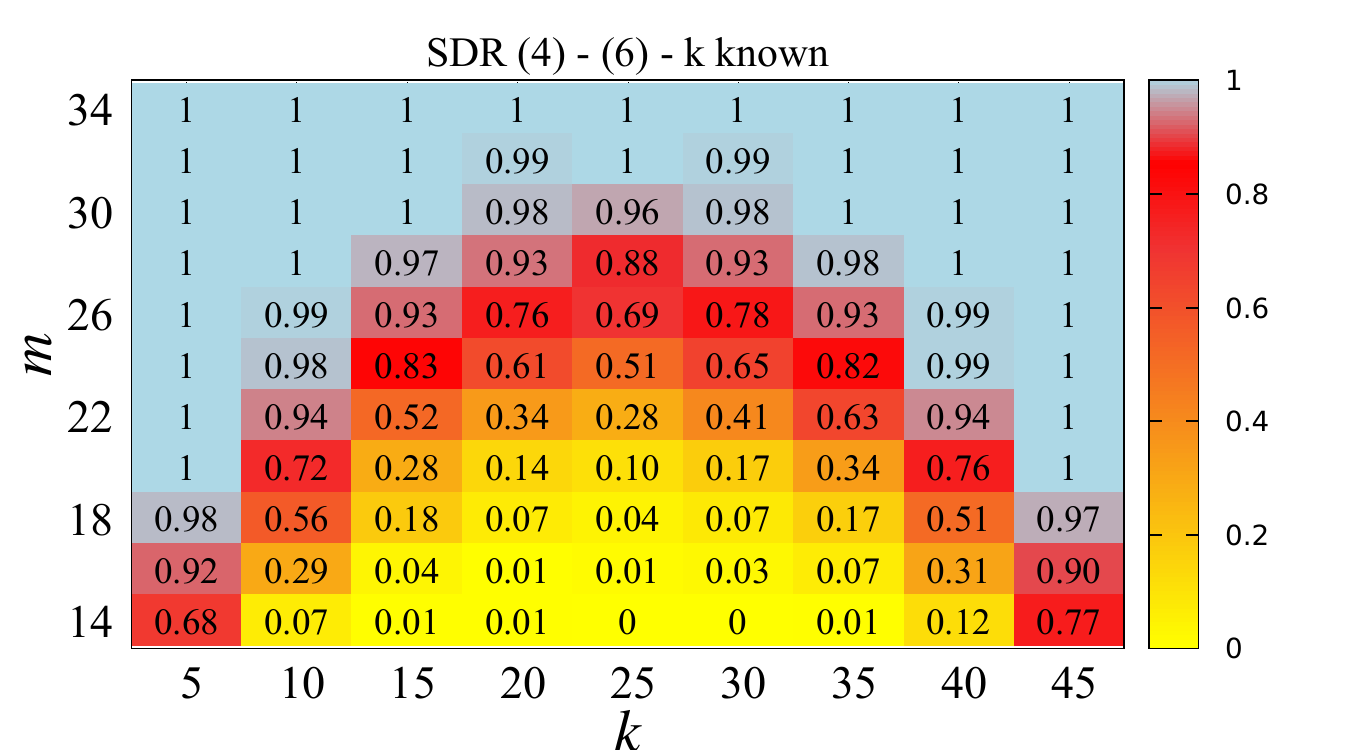}\\
\includegraphics[width=0.6\columnwidth]{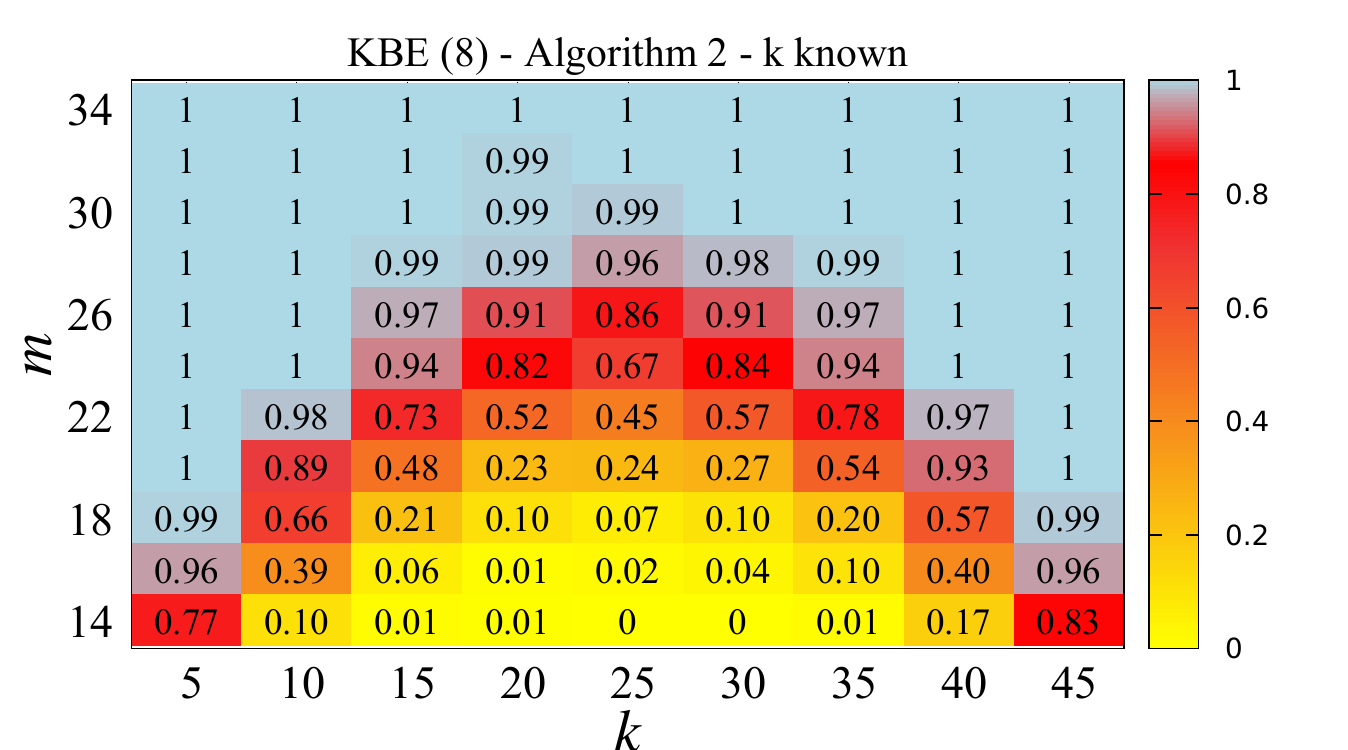}\\
\includegraphics[width=0.6\columnwidth]{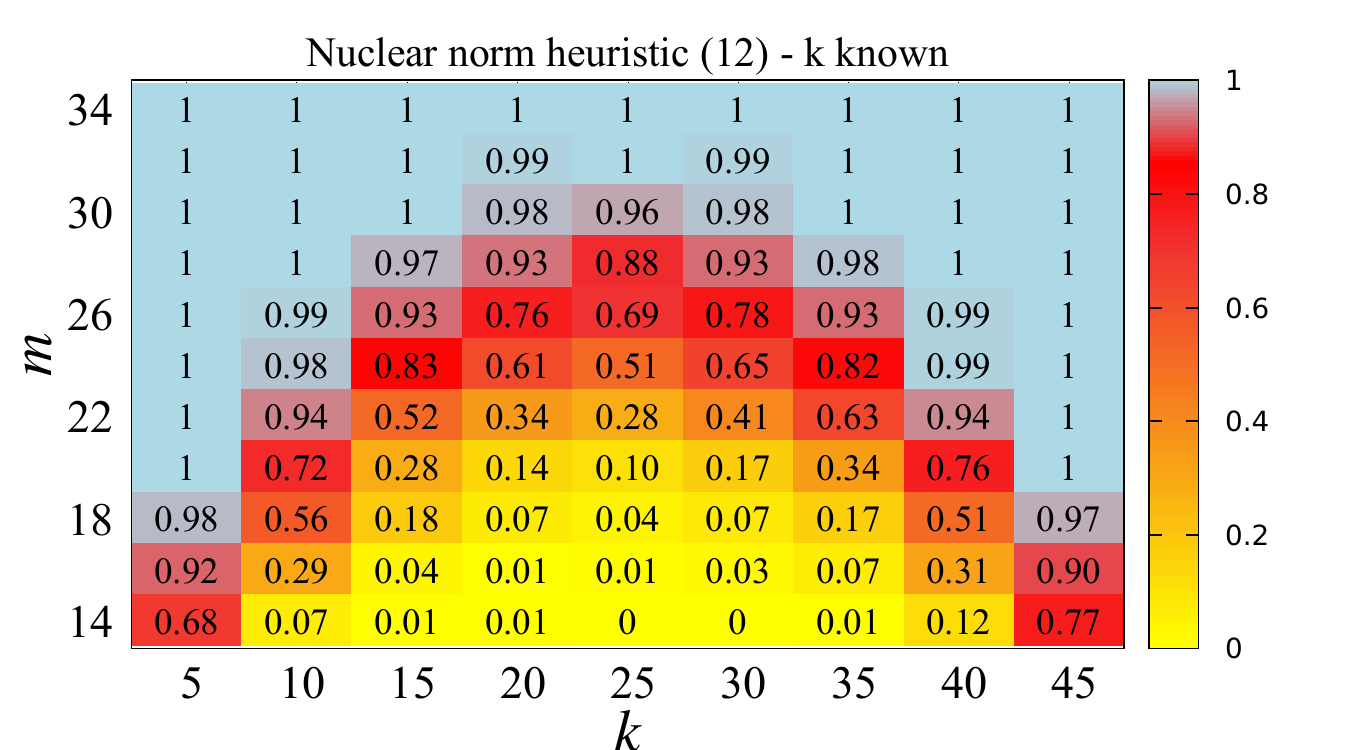}\\
\includegraphics[width=0.6\columnwidth]{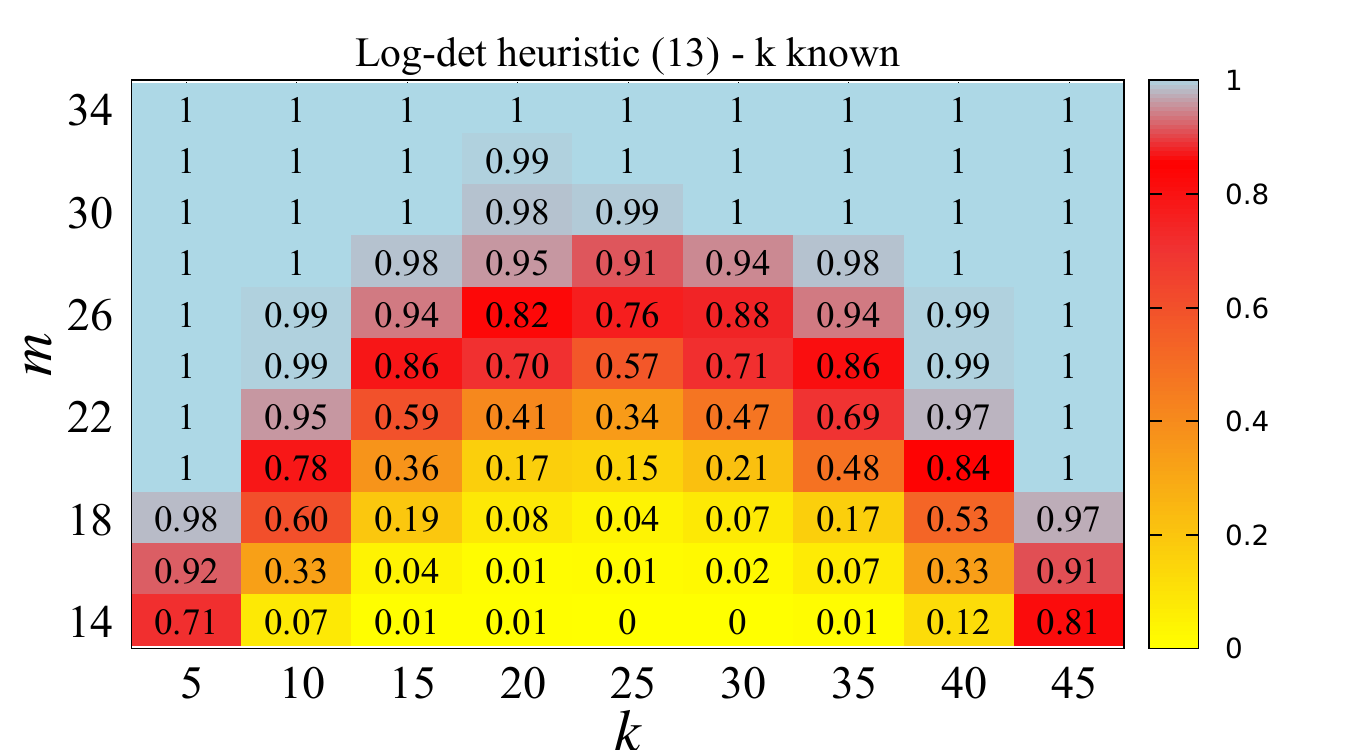}\\
\includegraphics[width=0.6\columnwidth]{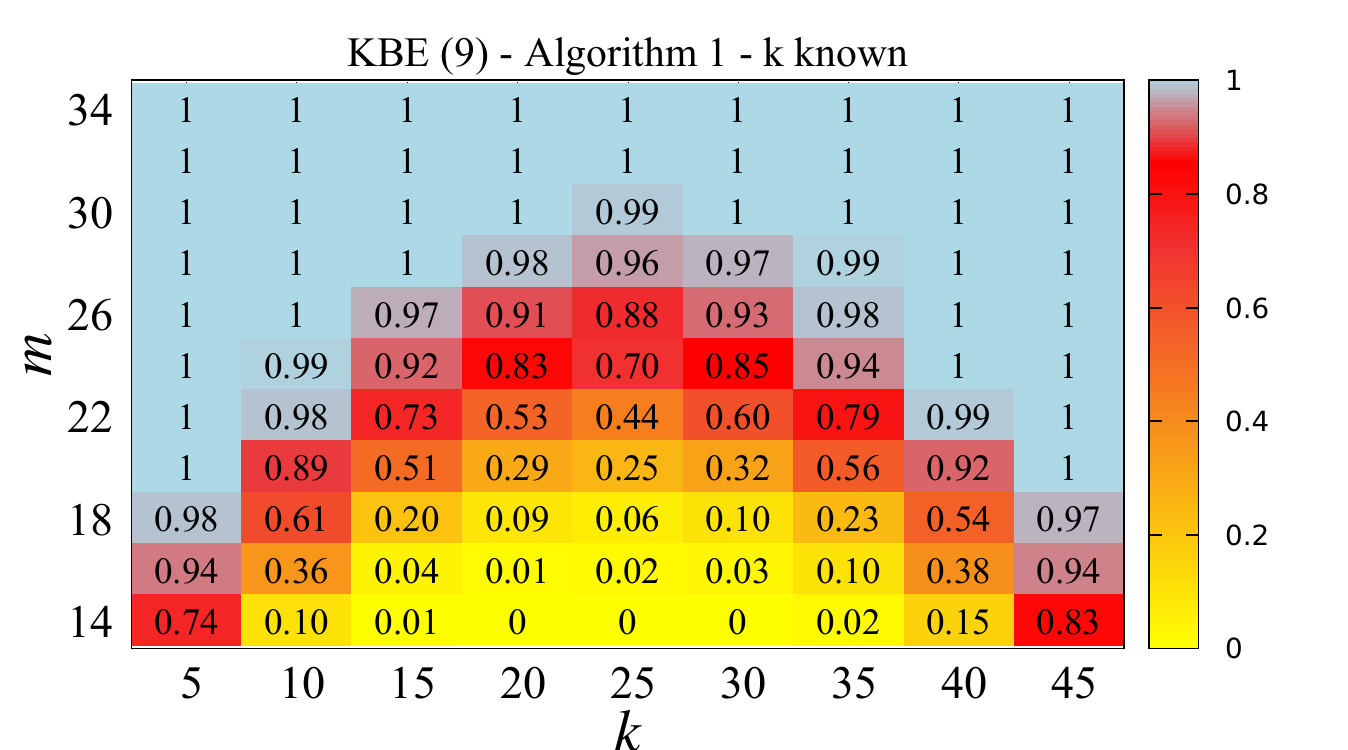}
\caption{Experiment 2 (known $k$): exact recovery rate.}
\label{fig:2}
\end{center}
\end{figure}
The considered setting is a follows. We take binary vectors $\xtrue\in\{0,1\}^n$, with $n=50$, sparsity level $k\in[5,45]$, and uniformly distributed support. $A$ is a random Gaussian matrix with $m$ rows, with $m\in[14,34]$. For all the methods, $\lambda=10^{-4}$, and  $T=3$. Algorithms \ref{alg:1} and \ref{alg:2}, and $\log$-$\det$ heuristic \eqref{sdr:ld} are initialized with the solutions of  SDR's \eqref{sdr:1} and \eqref{sdr:2}. If a binary solution is not achieved, i.e., the algorithm stops in a local minimum or saddle point, we randomly reinitialize for a maximum of 5 times. In \eqref{sdr:ld}, we set $\epsilon=10^{-6}$. To solve the involved SDR's, we use the Mosek C++ Fusion API, see \cite{mosek}, which guarantees fast solutions even for quite large dimensional problems. For reasons of space, a complete analysis of the complexity and large-scale experiments are left for future work. We just remark that the matrix inversion in the $\log$-$\det$ heuristic yields an higher computational complexity when compared to Algorithms \ref{alg:1} and \ref{alg:2}, and might be prohibitive in large-scale problems.
The considered algorithms are compared in terms of exact recovery rate, i.e., the number of experiments where a perfect recovery is achieved. Results are averaged over 200 random runs. As SDR's \eqref{sdr:1} and \eqref{sdr:2} are different formulations of the same problem, they are conveyed in the same  graph.

In Fig. \ref{fig:1}, we show the case of unknown $k$. We observe that Algorithm \ref{alg:2}  increases the success rate with respect to SDR's \eqref{sdr:1}-\eqref{sdr:2}. Specifically, KBE - Algorithm \ref{alg:2} has a transition phase from low to high recovery rate (namely, from 60$\%$ to 90$\%$, highlighted by red color in the figure) with less measurements than SDR's \eqref{sdr:1}-\eqref{sdr:2}: $m\in[24,26]$ instead of $m\in[26,28]$. We remark that these two methods have approximately constant performance in $k$. In contrast, as discussed above, the nuclear norm and $\log$-$\det$ heuristics as well as KBE - Algorithm \ref{alg:1}  penalize the trace, hence they have a sparsifying effect that makes them favorable when $k$ is small. Among them, the proposed KBE - Algorithm \ref{alg:1} generally achieves the best performance. In conclusion KBE - Algorithm \ref{alg:1} is the most reliable approach in non sparse problems, while KBE - Algorithm \ref{alg:2} is the most reliable approach for sparse problems. This attests that the KBE methodology is effective.

In Fig. \ref{fig:2}, the same experiment  is reported in  case of known $k$. As discussed in Sec. \ref{sec:PA}, Algorithm \ref{alg:1} can specifically exploit the information on $k$. Moreover,  the equation $\one^T x = k$ is added to $Ax=b$ for all the methods. In Fig. \ref{fig:2}, we observe that, as expected, all the methods take advantage of the knowledge of $k$, in particular when $k$ is either very small or  very large $k$. Also in this experiment, we observe that the proposed KBE Algorithms \ref{alg:1}-\ref{alg:2} are more accurate than the other methods. In particular, they are successful in at least $99\%$ of runs when $m\geq 30$, which is not achieved by the competitors.

Finally, we observe that, in our experiments, the run-time is between 0.1 and 0.8 seconds for KBE Algorithms \ref{alg:1}-\ref{alg:2}, and between 0.2 and 1.4 seconds for the $\log$-$\det$ heuristic.
\section{Conclusion}\label{sec:C}
In this paper, we propose a method to enhance the semidefinite programming relaxations of Boolean quadratic problems.  These relaxations provide the right minimizer if a 1-rank solution is found; nevertheless, enforcing low-rank solutions is NP-hard. The proposed strategy overcomes this drawback by leveraging known information on the eigenvalues. The global minimum of the proposed cost functional is proven to be the exact solution.  Since the problem is non-convex, a low-complex descent algorithm is developed, and shown to be effective through numerical results. The accuracy is enhanced with respect to state-of-the-art methods.
Future work will envisage a formal proof of the convergence of the proposed algorithms, and the development of different descent techniques, with particular focus on large scale problems.
%
%
\bibliographystyle{plain}
\bibliography{refs}             
\end{document}